\documentstyle[11pt,amssymb]{article}
\begin{document}
\noindent {\large \bf Harmonics on the quantum Euclidean space
related\\ to the quantum orthogonal group}
 \vskip 10 pt

N. Z. Iorgov and A. U. Klimyk

{\it Bogolyubov Institute for Theoretical Physics, Kiev 03143,
Ukraine}

\vskip  40 pt

\begin{abstract}
The aim of this paper is to study harmonic polynomials on the
quantum Euclidean space $E^N_q$ generated by elements $x_i$,
$i=1,2,\cdots ,N$, on which the quantum group $SO_q(N)$ acts. The
harmonic polynomials are defined as solutions of the equation
$\Delta _qp=0$, where $p$ is a polynomial in $x_i$, $i=1,2,\cdots
,N$, and the $q$-Laplace operator $\Delta _q$ is determined in
terms of the differential operators on $E^N_q$. The projector
${\sf H}_{m}: {\cal A}_{m}\to {\cal H}_{m}$ is constructed, where
${\cal A}_{m}$ and ${\cal H}_{m}$ are the spaces of homogeneous of
degree $m$ polynomials and homogeneous harmonic polynomials,
respectively. By using these projectors, a $q$-analogue of the
classical zonal polynomials and associated spherical polynomials
with respect to the quantum subgroup $SO_q(N-2)$ are constructed.
The associated spherical polynomials constitute an orthogonal
basis of ${\cal H}_{m}$. These polynomials are represented as
products of polynomials depending on $q$-radii and $x_j$,
$x_{j'}$, $j'=N-j+1$. This representation is in fact a
$q$-analogue of the classical separation of variables. The dual
pair $(U_q({\rm sl}_2), U_q({\rm so}_n))$ is related to the action
of $SO_q(N)$ on $E^N_q$. Decomposition into irreducible
constituents of the representation of the algebra $U_q({\rm
sl}_2)\times U_q({\rm so}_n)$ defined by the action of this
algebra on the space of all polynomials on $E^N_q$ is given.
\end{abstract}

\noindent {\sf I. INTRODUCTION}
\medskip

The Laplace operator, harmonic polynomials and related separations
of variables are of a great importance in classical analysis. They
are closely related to the rotation group $SO(N)$ and its
subgroups (see, for example, [1], chapter 10).

Harmonic polynomials are defined by the equation $\Delta p=0$,
where $\Delta$ is the Laplace operator and $p$ belongs to the
space ${\cal R}$ of polynomials on the Euclidean space $E_N \sim
{\Bbb R}^N$. The space ${\cal H}$ of all harmonic polynomials on
$E_N$ decomposes as a direct sum of the subspaces ${\cal H}_m$ of
homogeneous harmonic polynomials of degree $m$: ${\cal
H}=\bigoplus _{m=0}^\infty {\cal H}_m$. The Laplace operator
$\Delta$ on the Euclidean space $E_N$ commutes with the natural
action of the rotation group $SO(N)$ on this space. This means
that the subspaces ${\cal H}_m$ are invariant with respect to
$SO(N)$. The irreducible representation $T_m$ of the group $SO(N)$
with highest weight $(m,0,\cdots ,0)$ is realized on ${\cal H}_m$.

The Laplace operator $\Delta$ permits separation of variables on
the space ${\cal H}_m$. In other words, there exist different
coordinate systems (spherical, polyspherical) on $E_N$ and for
each of them it is possible to find the corresponding basis of the
space of solutions of the equation $\Delta p=0$ consisting of
products of functions depending on separated variables. To
different coordinate systems there correspond different
separations of variables. From the other side, to different
coordinate systems there correspond different chains of subgroups
of the group $SO(N)$ (see [1], chapter 10, for detail of this
correspondence). The basis of the space ${\cal H}_m$ in separated
variables (for a fixed coordinate system) consists of products of
Jacobi polynomials multiplied by $r^{2m}$, where $r$ is the
radius. These polynomials (considered only on the sphere
$S^{N-1}$) are matrix elements of the class 1 (with respect to the
subgroup $SO(N-1)$) irreducible representations $T_m$ of $SO(N)$
belonging to the zero column.

Many new directions of contemporary mathematical physics are
related to quantum groups and noncommutative geometry. It is
natural to generalize the theory described above to noncommutative
spaces. Such generalizations can be of a great importance for
further development of some branches of mathematical and
theoretical physics related to noncommutative geometry.

The aim of this paper is to construct a $q$-deformation of some
aspects of the classical theory described above. In the $q$-case,
instead of the Euclidean space we have the quantum Euclidean
space. It is defined in terms of the associative algebra ${\cal
A}$ generated by the elements $x_1,x_2,\cdots , x_N$ satisfying
the certain defining relations. These elements play a role of
Cartesian coordinates of $E_N$.

The $q$-Laplace operator $\Delta _q$ on ${\cal A}$ is defined in
terms of $q$-derivatives (see formula (17) below). Instead of the
group $SO(N)$ we have the quantum group $SO_q(N)$ or the
corresponding quantum algebra $U_q({\rm so}_N)$. In our exposition
it is more convenient to deal with the algebra $U_q({\rm so}_N)$.
The $q$-harmonic polynomials on the quantum Euclidean space are
defined as elements $p$ of ${\cal A}$ (that is, polynomials in
quantum coordinates $x_1,x_2,\cdots ,x_N$) for which $\Delta _q
p=0$. By using the algebra $U_q({\rm so}_N)$ or the quantum group
$SO_q(N)$ it is possible to construct for $q$-harmonic polynomials
the theory similar to the theory for classical harmonic
polynomials described above. Namely, we construct projectors ${\sf
H}_m: {\cal A}_m\to {\cal H}_m$, where ${\cal A}_m$ and ${\cal
H}_m$ are the subspaces of homogeneous (of degree $m$) polynomials
in ${\cal A}$ and in the space ${\cal H}$ of all $q$-harmonic
polynomials from ${\cal A}$, respectively. Using these projectors
we construct in ${\cal H}_m$ a $q$-analogue of associated
spherical harmonics with respect to the quantum subgroup
$SO_q(N-2)$. They constitute an orthogonal basis of the space
${\cal H}_m$ corresponding to the chain of the quantum subgroups
$$
SO_q(N)\supset SO_q(N-2)\supset SO_q(N-4)\supset \cdots \supset
SO_q(3)\ ({\rm or}\ SO_q(2)).
$$
Here we obtain a $q$-analogue of the corresponding spherical
separated coordinates. Our construction is similar to one used by
us in [2] for the case of quantum complex vector space with the
quantum unitary group $U_q(N)$ as a quantum motion group.

The operator $\Delta _q$ and the operator ${\hat Q}$ of
multiplication by the squared $q$-radius together with the certain
operator equivalent to operator of homogeneity degree on the sets
of homogeneous polynomials generate the quantum algebra $U_q({\rm
sl}_2)$. Thus, the algebra $U_q({\rm sl}_2)\times U_q({\rm so}_N)$
acts on the space ${\cal A}$. We decompose this representation
into irreducible ones. The pair $(U_q({\rm sl}_2), U_q({\rm
so}_N))$ constitutes a dual pair of quantum algebras which is a
$q$-analogue of the corresponding classical dual pair.

Our constructions use essentially the results of paper [3], where
the operator $\Delta _q$ and the spaces ${\cal H}_m$ were defined.

Everywhere below we suppose that $q$ is not a root of unity. We
shall use two different definitions of $q$-numbers:
$$
[a]=\frac{1-q^a}{1-q},\ \ \ \ [a]_q=\frac{q^a-q^{-a}}{q-q^{-1}}.
$$
It is necessary to pay attention which of these definitions is
used in each concrete case.
\bigskip

\noindent {\sf II. THE QUANTUM ALGEBRAS $U_q({\rm so}_N)$ AND THE
QUANTUM EUCLIDEAN SPACE}
\medskip

It is well known that the rotation group $SO(N)$ naturally acts on
the $N$-dimensional Euclidean space $E_N$. All this machinery (the
group $SO(N)$, the Euclidean space $E_N$, the action of $SO(N)$ on
$E_N$, etc) can be quantized. As a result, we have a "quantum"
action of the quantum rotation group $SO_q(N)$ on the quantum
Euclidean space (see [4]). The quantum Euclidean space $E^N_q$ is
defined by means of the algebra of polynomials ${\cal A}\equiv
{\Bbb C}_q[x_1, x_2,\cdots ,x_N]$ in noncommutative elements $x_1,
x_2,\cdots ,x_N$ which are called quantum Cartesian coordinates
(see [4] and [5]). The number $N$ can be even or odd and we
represent it as $N=2n$ or $N=2n+1$, respectively. Moreover, for
$j=1,2,\cdots ,N$ we shall use the notation $j'=N-j+1$.

The algebra ${\cal A}$ is the associative algebra generated by
elements $x_1,x_2,\cdots ,x_N$ satisfying the defining relations
$$
x_ix_j=qx_jx_i,\ \ \ \  i<j\ \ {\rm and}\ \  i\ne j', \eqno (1)
$$   $$
x_{i'}x_i-x_ix_{i'}=\frac{q-q^{-1}}{q^{\rho_i-1}+q^{-\rho_i+1}}
\sum _{j=i+1}^{(i+1)'} x_jx_{j'}q^{\rho_{j'}},\ \ \ i<n, \eqno (2)
$$   $$
x_{n'}x_n-x_nx_{n'}=(q^{1/2}-q^{-1/2})x^2_{n+1}\ \ {\rm if}\ \
N=2n+1 , \eqno (3)
$$   $$
x_{n'}x_n=x_nx_{n'} \ \ {\rm if}\ \ N=2n ,  \eqno (4)
$$
where
$$  {\textstyle
(\rho_1,\cdots ,\rho_{2n+1})=(n-\frac 12 ,n-\frac 32,\cdots ,\frac
12, 0,-\frac 12,\cdots ,-n+\frac 12) \ \ {\rm if}\ \ N=2n+1 ,}
$$  $$
(\rho_1,\cdots ,\rho_{2n})=(n-1 ,n-2,\cdots ,1,0, 0,-1,\cdots
,-n+1) \ \ {\rm if}\ \ N=2n .
$$
The set of monomials ${\bf x}^\nu:=x_1^{\nu_1}x_2^{\nu_2}\cdots
x_N^{\nu_N}$, $\nu_i=0,1,2\cdots$, form a basis of the algebra
${\cal A}$ (see [3]).

The vector space of the algebra ${\cal A}$ can be represented as a
direct sum of the vector subspaces ${\cal A}_m$ consisting of
homogeneous polynomials of homogeneity degree $m$,
$m=0,1,2,\cdots$:
$$
{\cal A}=\bigoplus _{m=0}^\infty {\cal A}_{m} .
$$

A $*$-operation (that is, an involutive algebra anti-automorphism)
can be defined on the algebra ${\cal A}$ turning it into a
$*$-algebra. This $*$-operation is uniquely determined by the
relations $x_i^*=q^{\rho_{i'}}x_{i'}$, $i=1,2,\cdots , N$.

The quantum rotation group $SO_q(N)$ and the corresponding
quantized universal enveloping algebra $U_q({\rm so}_N)$ act on
the algebra ${\cal A}$. These actions are determined by each
other. It will be convenient for us to use the action of the
algebra $U_q({\rm so}_N)$. The last algebra is the Hopf algebra
generated by the elements $K_i,K_i^{-1},E_i,F_i$, $i=1,2,\cdots,
n$, satisfying the certain defining relations (see, for example,
section 6.1.3 in [5]), where $n$ is an integral part of $N/2$. The
algebra $U_q({\rm so}_N)$ is supplied by the Hopf algebra
operations. We adopt these operations determined in [3]. The
action of $X\in U_q({\rm so}_N)$ on an element $a\in {\cal A}$
will be denoted as $X\vartriangleright a$.

A $*$-operation can also be introduced on $U_q({\rm so}_N)$ (see,
for example, [5]) which determines the compact form of $U_q({\rm
so}_N)$. We denote this compact form by $U_q({\rm so}(N))$. The
action of $U_q({\rm so}(N))$ on the $*$-algebra ${\cal A}$ is
compatible with the $*$-action, that is, $(X\triangleright
a)^*=S(X)^*\triangleright a^*$ for $X\in U_q({\rm so}(N))$ and
$a\in {\cal A}$, where $S$ is the antipode on $U_q({\rm so}_N)$
(see [3]).

The action of $U_q({\rm so}_N)$ on ${\cal A}$ is explicitly given
in [3], Lemma 2.5. For $U_q({\rm so}_{2n+1})$ and $U_q({\rm
so}_{2n})$, the action of elements $E_k$ and $F_k$, $k=1,2,\cdots
,n-1$, are determined as
$$
E_k\triangleright {\bf x}^\nu=[\nu_{k+1}]_qq^{\nu_k-\nu_{k+1}+1}
{\bf x}^{\nu+\varepsilon_k-\varepsilon_{k+1}} -[\nu_{k'}]_q
q^{\nu_{k}-\nu_{k+1}-\nu_{k'}+\nu_{(k+1)'} +1} {\bf
x}^{\nu+\varepsilon_{(k+1)'}-\varepsilon_{k'}} ,
$$  $$
F_k\triangleright {\bf x}^\nu=[\nu_{k}]_q
q^{-\nu_k+\nu_{k+1}-\nu_{(k+1)'}+\nu_{k'}+1} {\bf
x}^{\nu-\varepsilon_k+\varepsilon_{k+1}}- [\nu_{(k+1)'}]_q
q^{-\nu_{(k+1)'}+\nu_{k'}+1} {\bf
x}^{\nu-\varepsilon_{(k+1)'}+\varepsilon_{k'}}.
$$
The action of elements $E_n$ and $F_n$ are given by the formulas
$$
E_n\triangleright {\bf x}^\nu=[\nu_{n+1}]q^{\nu_n-\nu_{n+1}+3/2}
{\bf x}^{\nu+\varepsilon_n-\varepsilon_{n+1}} -[\nu_{n+2}]_q
q^{\nu_{n}-\nu_{n+2}+1} {\bf
x}^{\nu+\varepsilon_{n+1}-\varepsilon_{n+2}} ,
$$  $$
F_n\triangleright {\bf x}^\nu=[\nu_{n}]_q q^{-\nu_n+\nu_{n+2}+1/2}
{\bf x}^{\nu-\varepsilon_n+\varepsilon_{n+1}}-
[\nu_{n+1}]q^{-\nu_{n+1}+\nu_{n+2}+1} {\bf
x}^{\nu-\varepsilon_{n+1}+\varepsilon_{n+2}}
$$
if $N=2n+1$ and by the formulas
$$
E_n\triangleright {\bf
x}^\nu=[\nu_{n+1}]_qq^{\nu_{n-1}-\nu_{n+1}+1} {\bf
x}^{\nu+\varepsilon_{n-1}-\varepsilon_{n+1}}-
[\nu_{n+2}]_qq^{\nu_{n-1}+\nu_n-2\nu_{n+1}-\nu_{n+2}+1} {\bf
x}^{\nu+\varepsilon_{n}-\varepsilon_{n+2}} ,
$$  $$
F_n\triangleright {\bf x}^\nu=[\nu_{n-1}]_q q^{-\nu_{n-1}-2\nu_n
+\nu_{n+1}+\nu_{n+2}+1} {\bf
x}^{\nu-\varepsilon_{n-1}+\varepsilon_{n+1}}-
[\nu_{n}]_qq^{-\nu_{n}+\nu_{n+2}+1} {\bf
x}^{\nu-\varepsilon_{n}+\varepsilon_{n+2}}
$$
if $N=2n$. In these formulas $\varepsilon_i$ is the vector with
$i$th coordinate equal to 1 and all others equal to 0.

The monomials ${\bf x}^\nu$ are weight vectors with respect to the
action of $U_q({\rm so}_N)$ on ${\cal A}$. We represent weights
$\lambda$ in the well known orthogonal coordinate system, that is,
as $\lambda =\mu_1\varepsilon_1+\mu_2\varepsilon_2+\cdots
+\mu_n\varepsilon_n$ (in this system highest weights are given by
the coordinates $\mu_1,\mu_2,\cdots ,\mu_n$ such that $\mu_1\ge
\mu_2\ge \cdots$). The weight of the monomial ${\bf x}^\nu$ is
$$
\lambda=(\nu_1-\nu_{1'})\varepsilon_1+(\nu_2-\nu_{2'})\varepsilon_2
+\cdots +(\nu_n-\nu_{n'})\varepsilon_n.
$$
The action of the element $K_i$, $i<n$, on the monomial ${\bf
x}^\nu$ is given by the formula
$$
K_i \triangleright {\bf x}^\nu=q^{(\nu_i-\nu_{i'})
-(\nu_{i+1}-\nu_{(i+1)'})} {\bf x}^\nu.
$$
Moreover, in the algebra $U_q({\rm so}_N)$ there exist elements
${\hat K}_i$, $i=1,2,\cdots ,n$, such that
$$
{\hat K}_i \triangleright {\bf x}^\nu=q^{(\nu_i-\nu_{i'})} {\bf
x}^\nu.   \eqno (5)
$$

A differential calculus is developed on the quantum Euclidean
space which is determined by the $R$-matrix of the quantum algebra
$U_q({\rm so}_N)$. There exist different formulations for this
differential calculus. We adopt the definition of the differential
operators $\partial_i$, $i=1,2,\cdots ,N$, used in [3]. These
operators act on the monomials ${\bf x}^\nu$ as
$$
\partial_k\triangleright {\bf x}^\nu=[\nu_k]_qq^{\nu_{k+1}+\cdots
+\nu_{1'}}{\bf x}^{\nu-\varepsilon_k},\ \ \ k\le n, \eqno (6)
$$   $$
\partial_{n+1}\triangleright {\bf x}^\nu=[\nu_{n+1}]q^{\nu_{n'}+\cdots
+\nu_{1'}}{\bf x}^{\nu-\varepsilon_{n+1}},\ \ \ {\rm if}\ \ \
N=2n+1,
$$  $$
\partial_{k'}\triangleright {\bf x}^\nu=[\nu_{k'}]_qq^{\nu_{k'}+
\nu_{(k-1)'}\cdots +\nu_{1'}}{\bf x}^{\nu-\varepsilon_{k'}}+
$$   $$
\sum _{j=k+1}^n [\nu_j]_q [\nu_{j'}]_q (q-q^{-1})
q^{\rho_k-\rho_j} q^{d_{kj}} {\bf
x}^{\nu+\varepsilon_{k}-\varepsilon_j-\varepsilon_{j'}}+
$$  $$
+ [\nu_{n+1}-1] [\nu_{n+1}] \frac{q-q^{-1}}{1+q} q^{\rho_k+2}
q^{e_{k}} {\bf x}^{\nu+\varepsilon_{k}-2\varepsilon_{n+1}},\ \ \
k\le n,   \eqno (7)
$$
where $d_{kj}=\nu_k+\cdots +\nu_{j-1}+\nu_{(j-1)'}+\cdots
+\nu_{1'}$ and $e_{k}=(\nu_k+\cdots +\nu_{1'})-2\nu_{n+1}$. The
last summand in (7) must be omitted for $N=2n$. The operators
$\partial_i$, $i=1,2,\cdots ,N$, satisfy the relations
$$
\partial_i\partial_j=q^{-1}\partial_j\partial_i\ \ \ \ i<j,\ \ i\ne j',
\eqno (8)
$$  $$
\partial_{i'}\partial_i-\partial_i\partial_{i'}=-\frac{q-q^{-1}}
{q^{{\rho_i}-1}+q^{-{\rho_i}+1}}\sum _{k=i+1}^{(i+1)'}\partial_k
\partial_{k'}q^{\rho_k},\ \ \ i<n, \eqno (9)
$$   $$
\partial_{n'}\partial_n-\partial_n\partial_{n'}=-(q^{1/2}-q^{-1/2})
\partial^2_{n+1}\ \ \ {\rm if}\ \ \ N=2n+1, \eqno (10)
$$   $$
\partial_{n'}\partial_n=\partial_n\partial_{n'} \ \ \ {\rm if}\ \ \
N=2n.
$$

The operators $\partial_k$ and the operators ${\hat x}_i$ of left
multiplication by $x_i$ satisfy certain relations which can be
represented by means of the quantum $R$-matrix of the algebra
$U_q({\rm so}_N)$. These relations are given in [3]. We need the
following ones:
$$
\partial_k{\hat x}_k={\hat
x}_k\partial_kq^{\delta_{kk'}-1}-(q-q^{-1})\sum_{j<k} {\hat
x}_j\partial_j +(q-q^{-1})\sigma_k q^{2\rho_{k'}}{\hat
x}_{k'}\partial_{k'}+c,\ \ \ k=1,2,\cdots N, \eqno (11)
$$  $$
\partial_k{\hat x}_j={\hat
x}_j\partial_k+(q-q^{-1})\sigma_{kj}q^{\rho_{j'}-\rho_{k}} {\hat
x}_{k'}\partial_{j'} ,\ \ \ k\ne j,j', \eqno (12)
$$  $$
\partial_k{\hat x}_{k'}=q {\hat x}_{k'} \partial_{k},\ \ \ k\ne
k', \ \ \ \ c {\hat x}_{k}=q{\hat x}_{k}c,\ \ \  c
\partial_k=q^{-1}\partial_k c,  \eqno (13)
$$
where $\sigma_k=1$ if $k>k'$ and $\sigma_k=0$ otherwise,
$\sigma_{kj}=1$ if $k>j'$ and $\sigma_{kj}=0$ otherwise, $c$ is
the linear operator which acts on the monomials ${\bf x}^\nu$ as
$c\triangleright {\bf x}^\nu=q^{\nu_1+\cdots +\nu_{1'}}  {\bf
x}^\nu$.
\bigskip

\noindent {\sf III. SQUARED $q$-RADIUS AND $q$-LAPLACE OPERATOR}
\medskip

The element
$$
Q=\sum _{i=1}^N q^{\rho_{i'}}x_ix_{i'} \eqno (14)
$$
of the algebra ${\cal A}$ is called the {\it squared $q$-radius}
on the quantum Euclidean space. It is an important element in
${\cal A}$. It is shown in [3] that the center of ${\cal A}$ is
generated by $Q$.

Using relations between the elements $x_j$ it is shown that
$$
Q=(1+q^{N-2})\biggl( \sum _{i=1}^n q^{\rho_{i'}}x_ix_{i'}+
\frac{q}{q+1} x^2_{n+1}\biggr) \ \ \ {\rm if}\ \ \ N=2n+1,
$$
$$
Q=(1+q^{N-2}) \sum _{i=1}^n q^{\rho_{i'}}x_ix_{i'} \ \ \ {\rm if}\
\ \ N=2n.
$$
We shall also use the elements
$$
Q_j=\sum _{i=j}^{j'} q^{\rho_{i'}}x_ix_{i'},\ \ \ 1<j\le n , \eqno
(15)
$$
which are squared $q$-radii for the subalgebras ${\Bbb
C}_q[x_j,\cdots , x_{j'}]$. They satisfy the relations
$$
Q_jQ_k=Q_kQ_j,\ \ \ \   x_ix_{i'}=q^{\rho_{i}}\left(
\frac{Q_i}{1+q^{N-2i}}-\frac{Q_{i+1}}{1+q^{N-2i-2}}\right) ,\ \ \
1\le i\le n,
$$    $$
x_iQ_j=q^{2}Q_jx_i,\ \ \ x_{i'}Q_j=q^{-2}Q_jx_{i'}\ \    {\rm for}
\ \ i<j,\ \ \ \  x_iQ_j=Q_jx_i \ \ {\rm for} \ \ j\le i\le j' .
$$
It can be checked by direct computation that
$$
x^k_1x^k_{1'}={Q'}^k_1 ( Q'_{2}/Q'_1;q^{2})_k, \eqno (16)
$$
where $Q'_1=q^{-\rho_{1'}}Q_1/(1+q^{N-2})$,
$Q'_{2}=q^{-\rho_{1'}}Q_{2}/(1+q^{N-4})$ and
$$
(a;q)_s=(1-a)(1-aq)\cdots (1-aq^{s-1}).
$$

To the element (15) there corresponds the operator ${\hat Q}$ on
${\cal A}$ defined as
$$
{\hat Q}=\sum _{i=1}^N q^{\rho_{i'}} {\hat x}_i {\hat x}_{i'},
$$
where ${\hat x}_i$ is the operator of left multiplication by
$x_i$. It is clear that ${\hat Q} : {\cal A}_{m}\to {\cal
A}_{m+2}$.

We also consider on ${\cal A}$ the operator
$$
\Delta _q =\sum _{i=1}^N q^{\rho_i}\partial_i \partial _{i'}
\eqno (17)
$$
which is called the $q$-{\it Laplace operator} on the quantum
Euclidean space. We have $\Delta_q : {\cal A}_{m}\to {\cal
A}_{m-2}$.

The important property of the operators ${\hat Q}$ and $\Delta_q$
is that they commute with the action of the algebra $U_q({\rm
so}_N)$ on ${\cal A}$ (see [3]).

The operators ${\hat Q}$ and $\Delta_q$ satisfy the relations
$$
\Delta_q {\hat Q}^k-q^{2k} {\hat Q}^k\Delta_q ={\hat
Q}^{k-1}q^{-N+3} [2k][N+2k+2\gamma
-2]\frac{(1+q^{N-2})^2}{(1+q)^2} , \eqno (18)
$$     $$
\Delta_q ( Q^k)= Q^{k-1}q^{-N+3}[2k][N+2k-2]
\frac{(1+q^{N-2})^2}{(1+q)^2} \eqno (19)
$$
where $\gamma$ is the operator acting on the monomials ${\bf
x}^\nu$ as $\gamma {\bf x}^\nu =(\nu_1+\cdots +\nu_N){\bf x}^\nu$
(see [3]). We shall also use the following formula from [3]:
$$
\Delta_q ({\bf x}^\nu)=(1+q^{N-2})q^{\nu_1+\cdots
+\nu_{1'}-1}\times
$$   $$
\times \biggl( \sum _{j=1}^n  [\nu_j]_q[\nu_{j'}]_q q^{-\rho_j}q^d
{\bf x}^{\nu-\varepsilon_j-\varepsilon_{j'}} +
[\nu_{n+1}-1][\nu_{n+1}]\frac{q^e}{1{+}q}{\bf
x}^{\nu-2\varepsilon_{n+1}} \biggr) , \eqno (20)
$$
where $d=\nu_1+\cdots +\nu_{j-1}+\nu_{(j-1)'}+\cdots +\nu_{1'}$,
$e=\nu_1+\cdots +\nu_{1'}-2\nu_{n+1}+2$, and the last summand must
be omitted for $N=2n$.
\bigskip

\noindent {\sf IV. $q$-HARMONIC POLYNOMIALS}
\medskip

A polynomial $p\in {\cal A}$ is called $q$-{\it harmonic} if
$\Delta_q p=0$. The linear subspace of ${\cal A}$ consisting of
all $q$-harmonic polynomials is denoted by ${\cal H}$. If ${\cal
H}_{m}={\cal A}_{m}\cap {\cal H}$, then ${\cal H}_{m}$ is the
subspace of ${\cal H}$ consisting of all homogeneous of degree $m$
harmonic polynomials.
\medskip

{\bf Proposition 1} [3]. {\it If $m\ge 2$, then the space ${\cal
A}_{m}$ can be represented as the direct sum}
$$
{\cal A}_{m}={\cal H}_{m}\oplus Q{\cal A}_{m-2}. \eqno (21)
$$

We shall need the following consequences of the decomposition
(21):
\medskip

{\it Corollary 1.} {\it If $p\in {\cal H}_{m}$, then $p$ cannot be
represented as $p=Q^kp'$, $k\ne 0$, with some polynomial $p'\in
{\cal A}$.}
\medskip

{\it Corollary 2.} {\it The space ${\cal A}_{m}$ decomposes into
the direct sum ${\cal A}_{m}=\bigoplus _{j=0}^{\lfloor m/2\rfloor}
Q^j{\cal H}_{m-2j}$, where $\lfloor m/2\rfloor$ is the integral
part of the number $m/2$.}
\medskip

{\it Corollary 3.} {\it For dimension of the space of $q$-harmonic
polynomials ${\cal H}_{m}$ we have the formula}
$$
{\rm dim}\ {\cal H}_{m}={(m+N-3)!(2m+N-2)\over (N-2)!m!} .
$$

{\it Corollary 4.} {\it The linear space ${\cal H}$ can be
represented in the form of a direct sum}
$$
{\cal H}=\bigoplus _{m=0}^\infty {\cal H}_{m} .
$$

Corollary 1 is a direct consequence of formula (21). Corollary 2
easily follows from repeated application of (21). Corollary 3 is
proved in the same way as in the classical case (see, for example,
[1], Chap. 10). For this we note that
$$
{\dim}\ {\cal A}_{m}={(N+m-1)!\over (N-1)!m!} .
$$
Hence, for ${\dim}\ {\cal H}_{m}={\dim}\ {\cal A}_{m}- {\dim}\
{\cal A}_{m-2}$ we obtain the expression stated in the corollary.
In order to prove Corollary 4 we note that any $p\in {\cal H}$ can
be represented as $p=\sum _m p_m$, $p_m\in {\cal A}_m$. We have
$\Delta_qp=\sum _m \Delta_q p_m=0$. Since $\Delta_q p_m$,
$m=0,1,2,\cdots$, have different homogeneity degrees, it follows
from the last equality that $\Delta_q p_m=0$ for all values of
$m$. Thus, ${\cal H}=\bigoplus _{m=0}^\infty {\cal H}_{m}$.
\medskip

{\it Remark:} If $n=2$, then ${\cal A}$ consists of all
polynomials in commuting elements $x_1$ and $x_{1'}\equiv x_2$. In
this case, the space ${\cal H}$ of $q$-harmonic polynomials has a
basis consisting of the polynomials
$$
1, \ \ x_1^k,\ \ x_{1'}^k,\ \ \ \ k=1,2,\cdots . \eqno (22)
$$

{\bf Proposition 2.} {\it The linear space isomorphism ${\cal
A}\simeq {\Bbb C}[Q]\otimes {\cal H}$ is true, where ${\Bbb C}[Q]$
is the space of all polynomials in $Q$.}
\medskip

This proposition follows from Corollary 2.
\medskip

The decomposition ${\cal A}\simeq {\Bbb C}[Q]\otimes {\cal H}$ is
a $q$-analogue of the theorem on separation of variables for Lie
groups in an abstract form (see [6]). It follows from this
decomposition that
$$
{\cal A}\simeq {\Bbb C}[Q]\otimes {\cal H}\simeq {\Bbb
C}[Q]\otimes \bigoplus _{m\ge 0}{\cal H}_{m} = \bigoplus _{m\ge 0}
\left( {\Bbb C}[Q]\otimes {\cal H}_{m} \right) . \eqno (23)
$$

Since the operator $\Delta_q$ commutes with the action of the
algebra $U_q({\rm so}_n)$ the subspaces ${\cal H}_{m}$ are
invariant with respect to the action of this algebra. It is proved
in [3] that the irreducible representation $T_m$ of $U_q({\rm
so}_N)$ with highest weight $(m,0,\cdots ,0)$ is realized on
${\cal H}_{m}$.

We denote by ${\cal A}^{U_q({\rm so}_r)}$ the space of elements of
${\cal A}$ consisting of invariant elements with respect to the
action of $U_q({\rm so}_r)$. We have ${\cal A}^{U_q({\rm
so}_N)}={\Bbb C}[Q]$ (see [3]). In what follows we shall consider
the subalgebra $U_q({\rm so}_{N-2})$ generated by the elements
$H_i,E_i,F_i$, $i=2,3,\cdots ,n$.
\medskip

{\bf Proposition 3.} {\it We have}
$$
{\cal A}^{U_q({\rm so}_{N-2})}\simeq \bigoplus_{k,l} {\Bbb
C}[Q_{2}] x_1^k x_{1'}^l \simeq \bigoplus_{k,l} {\Bbb C}[Q] x_1^k
x_{1'}^l.
$$

{\it Proof.} In order to prove this proposition we note that for
$U_q({\rm so}_{N-2})$-module ${\cal A}$ we have
$$
{\cal A}={\Bbb C}_q[x_1,x_2,\cdots ,x_N]=\bigoplus _{k,l} {\Bbb
C}_q[x_2,x_3,\cdots ,x_{N-1}] x^k_1x^l_N .
$$
The action of $U_q({\rm so}_{N-2})$ on monomials $x^k_1x^l_N$ is
trivial. Moreover, ${\Bbb C}_q[x_2,x_3,\cdots ,x_{N-1}]^{U_q({\rm
so}_{N-2})} ={\Bbb C}[Q_{2}]$. Since $Q=c_1Q_{2}+c_2x_1x_N$, where
$c_1$ and $c_2$ are constants, we have ${\cal A}^{U_q({\rm
so}_{N-2})}\simeq \bigoplus_{k,l} {\Bbb C}[Q_{2}] x_1^k x_N^l
\simeq \bigoplus_{k,l} {\Bbb C}[Q] x_1^k x_N^l$. Proposition is
proved.
\bigskip

\noindent {\sf V. THE DUAL PAIR $(U_q({\rm sl}_2),U_q({\rm
so}_N))$}
\medskip

The formulas
$$
ke=q^2ek,\ \ \ \ kf=q^{-2}fk,\ \ \ \
ef-fe=\frac{k-k^{-1}}{q-q^{-1}} \eqno (24)
$$
determine the quantum algebra $U_q({\rm sl}_2)$ generated by the
elements $k,k^{-1},e,f$. Let ${\cal L}({\cal A})$ be the space of
linear operators on the algebra ${\cal A}$. It was shown in [3]
that the operators
$$
\omega (k)=q^{N/2}q^\gamma, \ \ \ \ \omega (e)={\hat Q} ,\ \ \ \
\omega (f)=-\Delta_q q^{-\gamma} \frac{q^{N/2}}{(1+q^{N-2})^2}
\eqno (25)
$$
satisfy relations (24). This means that the algebra homomorphism
$\omega : U_q({\rm sl}_2)\to {\cal L}({\cal A})$, uniquely
determined by formulas (25), is a representation of $U_q({\rm
sl}_2)$.

Since the operators $\omega (k)$, $\omega (e)$, $\omega (f)$
commute with the operators $L(X)$, $X\in U_q({\rm so}_N)$, we can
introduce the representation $\omega \times L$ of the algebra
$U_q({\rm sl}_2)\times U_q({\rm so}_N)$ on ${\cal A}$, where $L$
is the above defined natural action of $U_q({\rm so}_N)$ on ${\cal
A}$. This representation is reducible. Let us decompose it into
irreducible constituents.

By (23), we have ${\cal A}=\bigoplus _{m\ge 0} ({\Bbb C}[Q]\otimes
{\cal H}_{m})$. The subspaces ${\Bbb C}[Q]\otimes {\cal H}_{m}$
are invariant under $U_{q}({\rm sl}_2)\times U_q({\rm so}_N)$,
since the space ${\Bbb C}[Q]$ is elementwise invariant under
$U_q({\rm so}_N)$, and for $f\in {\Bbb C}[Q]$ and $h_{m}\in {\cal
H}_{m}$ we have
$$
\omega (e) (f(Q)\otimes h_{m})=Qf(Q)\otimes h_{m}, \eqno (26)
$$      $$
\omega (f)(Q^r\otimes h_{m})=-[r]_q[r+m+(N/2)-1]_q Q^{r-1}\otimes
h_{m}, \eqno (27)
$$    $$
\omega (k) (Q^{r}\otimes h_{m})=q^{2r+m+N/2}(Q^{r}\otimes h_{m})
\eqno (28)
$$
(we used formula (18) for obtaining (27)). These formulas show
that $U_{q}({\rm sl}_2)$ acts on ${\Bbb C}[Q]$ and $U_q({\rm
so}_N)$ acts on ${\cal H}_{m}$. However, this action of
$U_{q}({\rm sl}_2)$ depends on the component ${\cal H}_{m}$.
Taking the basis
$$
| r\rangle :=[r+m+(N/2)-1]_q!^{-1}Q^r, \ \ \ \ r=0,1,2,\cdots ,
$$
in the space ${\Bbb C}[Q]$, we find from (26)--(28) that
$$
\omega (k)|r\rangle =q^{2r+m+N/2}|r\rangle \ \ \ \ \ \omega (f)|
r\rangle =-[r]_q| r-1\rangle ,
$$    $$
\omega (e) |r\rangle =[r+m+N/2]_q|r+1\rangle .
$$
Comparing this representation with the known irreducible
representations of $U_{q}({\rm sl}_2)$ (see, for example, [7]) we
derive that the irreducible representation of $U_{q}({\rm sl}_2)$
of the discrete series with lowest weight $m+N/2$ is realized on
the component ${\Bbb C}[Q]$ of the space ${\Bbb C}[Q]\otimes {\cal
H}_{m}$. We denote this representation of $U_{q}({\rm sl}_2)$ by
$D_{m+N/2}$.

Thus, we have derived that on the subspace ${\Bbb C}[Q]\otimes
{\cal H}_{m}\subset {\cal A}$ the irreducible representation
$D_{m+N/2}\times T_{m}$ of the algebra $U_{q}({\rm sl}_2)\times
U_q({\rm so}_N)$ acts. This means that for the reducible
representation $\omega \otimes L$ we have the following
decomposition into irreducible components:
$$
\omega \times L=\bigoplus ^\infty _{m=0} D_{m+N/2}\times T_{m},
$$
that is, each irreducible representation of $U_q({\rm so}_N)$ in
this decomposition determines uniquely the corresponding
irreducible representation of $U_q({\rm sl}_2)$ and vise versa.
This means that $U_q({\rm sl}_2)$ and $U_q({\rm so}_N)$ constitute
a {\it dual pair} under the action on ${\cal A}$. It is a
$q$-analogue of the well known classical dual pair $({\rm sl}_2,
{\rm so}_N)$ (see, for example, [1], Chapter 12).
 \bigskip

\noindent {\sf VI. RESTRICTION OF $q$-HARMONIC POLYNOMIALS ONTO
THE QUANTUM SPHERE}
 \medskip

The associative algebra ${\cal F}(S_q^{N-1})$ generated by the
elements $x_1,\cdots ,x_N$ satisfying the relations (1)--(3) and
the relation $Q=1$ is called {\it the algebra of functions on the
quantum sphere} $S_q^{N-1}$ (see [4] and [5], Chap. 11). It is
clear that the following canonical algebra isomorphism has place:
$$
{\cal F}(S_q^{N-1}) \simeq {\cal A}/{\cal I} ,
$$
where ${\cal I}$ is the two-sided ideal of ${\cal A}$ generated by
the element $Q-1$. We denote by $\tau$ the canonical algebra
homomorphism
$$
\tau : {\cal A}\to {\cal A}/{\cal I}\simeq {\cal F}(S_q^{N-1}) .
$$
This homomorphism is called the {\it restriction} of polynomials
of ${\cal A}$ onto the quantum sphere $S_q^{N-1}$.

It was shown in [3] that $\tau :{\cal H}\to {\cal F}(S_q^{N-1})$
is a one-to-one mapping, that is, the restriction of a
$q$-harmonic polynomial to the sphere $S_q^{N-1}$ determines this
polynomial uniquely. This statement allows us to determine a
scalar product on ${\cal H}$. For this, we use the invariant
functional $h$ on the quantum sphere $S_q^{N-1}$ defined in [3].
In order to give this functional we introduce the linear subspace
$(\tau {\cal A})^0$ of ${\cal F}(S_q^{N-1})$ spanned by the
elements $\tau {\bf x}^\nu$ such that
$$
\nu_1=\nu_{1'},\cdots , \nu_n=\nu_{n'},\ \  \nu_{n+1}=2m,\ \
m=0,1,2,\cdots,
$$
(for $N=2n$ the last condition must be omitted). The functional
$h$ vanishes on the elements $\tau {\bf x}^\nu \not\in (\tau {\cal
A})^0$ and on the monomials $\tau {\bf x}^\nu\in (\tau {\cal
A})^0$ it is given by the formula
$$
h(\tau {\bf x}^\nu)=\frac{(q^{-2};q^{-2})_{\nu _1}\cdots
(q^{-2};q^{-2})_{\nu _n} (q^{-1};q^{-2})_m(1+q)^m}
{q^{(\rho_1\nu_1+\cdots +\rho_n\nu_n)+m}(1+q^{N-2})^{\nu_1+\cdots
+\nu_n+m}(q^{-N};q^{-2})_{\nu _1+\cdots +\nu_n+m}},
$$
where, as before, $m=0$ for $N=2n$. The following assertions
(similar to ones proved in [8] for the case of the quantum group
$U_q(N)$) are true:

(a) The subalgebra $(\tau {\cal A})^0$ is a commutative algebra
generated by $Q_{n-1}$, $Q_{n-2},\cdots , Q_1$ and also by
$x_{n+1}^2$ if $N=2n+1$.

(b) The algebra $(\tau {\cal A})^0$ is isomorphic to the
polynomial algebra in $n-1$ commuting indeterminates if $N=2n$ and
in $n$ commuting indeterminates if $N=2n+1$ .

A scalar product $\langle \cdot ,\cdot \rangle$ on ${\cal H}$ is
introduced by the formula:
$$
\langle p_1,p_2\rangle =h((\tau p_1)^*(\tau p_2)), \eqno (29)
$$
where $a^*$ determines an element conjugate to $a\in {\cal A}$
under action of the $*$-operation introduced in section 2.
\medskip

{\bf Proposition 4.} {\it We have ${\cal H}_{m}\bot {\cal H}_{r}$
if $m\ne r$.}
\medskip

{\it Proof} follows from the fact that $(\tau p_1)^*(\tau
p_2)\not\in (\tau {\cal A})^0$ if $p_1\in {\cal H}_{m}$, $p_2\in
{\cal H}_{r}$, and $m\ne r$.
\bigskip

\noindent {\sf VII. THE PROJECTION ${\cal A}_{m}\to {\cal H}_{m}$}
 \medskip

Let us go back to the decomposition (21) and construct the
projector
$$
{\sf H}_{m} :{\cal A}_{m} ={\cal H}_{m}\oplus Q{\cal A}_{m-2}\to
{\cal H}_{m}.
$$
We present this projector in the form
$$
{\sf H}_{m} p=\sum _{k=0}^{\lfloor m/2\rfloor } \alpha _k {\hat
Q}^k\Delta_q^k p, \ \ \ \  \alpha _k \in {\Bbb C},\ \ \ p\in {\cal
A}_{m}, \eqno (30)
$$
where $\lfloor m/2\rfloor$ is the integral part of the number
$m/2$. Let us show that the summands on the right hand side are
linearly independent at least for one nontrivial $p\in {\cal
A}_{m}$ (in this case the coefficients $\alpha_k$ are determined
uniquely up to a common constant). Let $p=x_{n+1}^m$ if $N=2n+1$.
Using formula (20) we derive that
$$
\Delta_q^k(x^m_{n+1})=q^{k}\left(  \frac{1+q^{N-2}}{1+q}\right)^k
\frac{[m]!}{[m-2k]!} x^{m-2k}_{n+1} ,\ \ \ \ 2k\le m, \eqno (31)
$$
where $[m]!=[1][2]\cdots [m]$. Then the right hand side of (30) is
a linear combination of the elements $Q^kx^{m-2k}_{n+1}$,
$k=1,2,\cdots, \lfloor m/2\rfloor$. It is easy to see that these
elements are linearly independent. If $N=2n$, then instead of
$p=x^m_{n+1}$ we take $p=x_1^{m_1}x_{1'}^{m'_1}$ and make the same
reasoning (see these calculation in the next section).

We have to calculate values of the coefficients $\alpha _k$ in
(30). In order to do this, we act by the operator $\Delta_q$ upon
both sides of (30) and use the relation (18). Under this action,
the left hand side vanishes. Equating the right hand side to 0 and
taking into account that the elements ${\hat Q}^k\Delta_q^{k+1}
p$, $k=1,2,\cdots,\lfloor m/2\rfloor$, are linearly independent
for generic elements $p\in {\cal A}_m$ , we derive the recurrence
relation
$$
q^{-N-2k+5}\frac{(1+q^{N-2})^2}{(1+q)^2} [2k][N+2m-2k-2] \alpha
_k+\alpha _{k-1}=0
$$
for $\alpha _k$ which gives
$$
\alpha _k =(-1)^kq^{(N-4)k+k^2}\frac{(1+q)^{2k}}{(1+q^{N-2})^{2k}}
\frac{[N+2m-2k-4]!!}{[2k]!![N+2m-4]!!},
$$
where $[s]!!=[s][s-2][s-4]\cdots [2]\ ({\rm or}\ [1])$ for $s\ne
0$ and $[0]!!=1$. Using the relations
$$
[2k]!!=\frac{(q^2;q^2)_k}{(1-q)^k},\ \ \ \
\frac{[N+2m-2k-4]!!}{[N+2m-4]!!}=\frac{(1-q)^k}{(q^{N+2m-2k-2};q^2)_k},
$$    $$
(q^{N+2m-2k-2};q^2)_k=
(q^{-N-2m+4};q^2)_kq^{-2k-k(k-1)}(-q^{N+2m-2})^k,
$$
we derive that
$$
\alpha_k=\frac{q^{2k^2-2mk-k}(1-q^2)^{2k}}{(1+q^{N-2})^{2k}
(q^{-N-2m+4};q^2)_k(q^2;q^2)_k}.   \eqno (32)
$$

Note that the coefficients $\alpha _k$ are determined by the
recurrence relation uniquely up to a common constant. In (30) we
have chosen this constant in such a way that ${\sf H}_{m}p=p$ for
$p\in {\cal H}_{m}$. This means that ${\sf H}_{m}^2={\sf H}_{m}$.
\medskip

{\bf Proposition 5.} {\it The operator ${\sf H}_{m}$ commutes with
the action of $U_q({\rm so}_N)$.}
\medskip

{\it Proof.} This assertion follows from the fact that the action
of $X\in U_q({\rm so}_N)$ commutes with ${\hat Q}$ and $\Delta_q$.
Proposition is proved.
\medskip

The operator ${\sf H}_{m}$ can be used for obtaining explicit
forms of $q$-harmonic polynomials. As an example, we derive here
formulas for harmonic projection of the polynomial $x_{n+1}^m\in
{\cal A}_m$ when $N=2n+1$. For this we use formula (31) for
$\Delta_q^k(x^m_{n+1})$. Since
$$
\frac{[m]!}{[m-2k]!}=\frac{(q^{m-2k+2};q^2)_k (q^{m-2k+1};q^2)_k}{
(1-q)^{-2k}}.
$$   $$
(q^{m-2k+2};q^2)_k =(-1)^k q^{mk-k(k-1)} (q^{-m};q^2)_k ,
$$   $$
(q^{m-2k+1};q^2)_k =(-1)^k q^{mk-k^2} (q^{-m+1};q^2)_k
$$
we derive that
$$
\Delta_q^k(x^m_{n+1})=\left( \frac{1+q^{N-2}}{1+q}\right)^k
\frac{q^{2mk-2k^2+2k}}{(1-q)^{2k}} (q^{-m};q^2)_k (q^{-m+1};q^2)_k
x^{m-2k}_{n+1}.
$$

Using the expression (30) for ${\sf H}_{m}x^m_{n+1}$ and formula
(32) for coefficients $\alpha_k$ we obtain
$$
{\sf H}_{m}x^m_{n+1}=x^m_{n+1}\sum _{k=0}^{\lfloor m/2\rfloor }
\frac{(q^{-m};q^2)_k (q^{-m+1};q^2)_k}{(q^2;q^2)_k
(q^{-N-2m+4};q^2)_k} (aQx_{n+1}^{-2})^k,   \eqno (33)
$$
where
$$
a=\frac{q(1+q)}{1+q^{N-2}}.
$$
Note that we used $x_{n+1}^{-2}$ in (33). However, since there
exists the multiplier $x^m_{n+1}$ before the sign of sum, negative
powers of $x_{n+1}$ in fact are absent.

The expression (33) for ${\sf H}_{m}x^m_{n+1}$ can be represented
in terms of the basic hypergeometric function ${}_2 \phi_1$ (see
[9] for the definition of this function):
$$
{\sf H}_{m}x^m_{n+1}=x^m_{n+1} {}_2 \phi_1 (q^{-m}, q^{-m+1};
q^{-N-2m+4}:\ q^2,\ aQx_{n+1}^{-2}).
$$
Using the definition
$$
P^{(\alpha ,\beta )}_k(x;\ q)= {}_2 \phi _1 (q^{-k},\ q^{\alpha
+\beta +k+1};\ q^{\alpha +1};\ q,\ qx)
$$
of the little $q$-Jacobi polynomials we can represent ${\sf
H}_{m}x^m_{n+1}$ in the form
$$
{\sf H}_{m}x^m_{n+1}=x^m_{n+1}P^{(-\frac N2
+m+1,\frac{N-3}{2})}_{m/2} \left( \frac{1+q}{q(1+q^{N-2})}
Qx^{-2}_{n+1}\right)
$$
if $m$ is even and in the form
$$
{\sf H}_{m}x^m_{n+1}=x^m_{n+1} P^{(-\frac N2
+m+1,\frac{N-3}{2})}_{(m-1)/2} \left( \frac{1+q}{q(1+q^{N-2})}
Qx^{-2}_{n+1}\right)
$$
if $m$ is odd.

Up to a constant the restriction of the expression for ${\sf
H}_{m}x^m_{n+1}$ to $S^{N-1}_q$ was found (by other method: as a
zonal spherical function with respect to a certain ideal) in [3].
It was expressed in term of the big $q$-Jacobi polynomials which
in fact are the basic hypergeometric functions ${}_3\varphi_2$ of
the argument $q$. Thus, we found another expression for these
zonal spherical functions.
\bigskip

\noindent {\sf VIII. ZONAL POLYNOMIALS WITH RESPECT TO
$SO_q(N-2)$}
 \medskip

A polynomial $\varphi$ of the space ${\cal H}_{m}$ is called {\it
zonal} with respect to the quantum subgroup $SO_q(N-2)$ (or with
respect to the subalgebra $U_q({\rm so}_{N-2})$) if it is
invariant with respect to action of elements $X\in U_q({\rm
so}_{N-2})$. In order to find zonal polynomials $\varphi \in {\cal
H}_{m}$ we have to take polynomials $p \in {\cal A}_{m}$ invariant
with respect to the subalgebra $U_q({\rm so}_{N-2})$ and to act on
them by the projection ${\sf H}_{m}$.

It follows from Proposition 3 that in the space ${\cal A}_{m}$
there exist $m+1$ elements which are $U_q({\rm
so}_{N-2})$-invariant and linearly independent over ${\Bbb C}[Q]$.
They coincide with $x_1^{m_1}x_{1'}^{m'_1}$, $m_1+m'_1=m$.
Therefore, ${\sf H}_m (x_1^{m_1}x_{1'}^{m'_1})$, $m_1+m'_1=m$, are
zonal polynomials with respect to $U_q({\rm so}_{N-2})$. Let us
find explicit form of these polynomials.

Using formula (20) we find that
$$
\Delta^k_q(x_1^{m_1}x_{1'}^{m'_1})=(1+q^{N-2})^kq^{(m-k)k}
q^{-(n-\epsilon)k} \frac{[m_1]_q![m'_1]_q!}{[m_1-k]_q![m'_1-k]_q!}
x_1^{m_1-k}x_{1'}^{m'_1-k},  \eqno (34)
$$
where $\epsilon=1$ for $N=2n$ and $\epsilon =\frac 12$ for
$N=2n+1$. Since
$$
\frac{[m_1]_q!}{[m_1-k]_q!}=q^{(2m_1-k+1)k/2}
\frac{(q^{-2m_1};q^2)_k}{(q-q^{-1})^k},
$$
we have
$$
\Delta^k_q(x_1^{m_1}x_{1'}^{m'_1})=(1+q^{N-2})^kq^{2(m-k)k}
q^{(-n+3+\epsilon)k} \frac{(q^{-2m_1};q^2)_k(q^{-2m'_1};q^2)_k}
{(1-q^2)^{2k}} x_1^{m_1-k}x_{1'}^{m'_1-k} .
$$
Now using formulas (30) and (32) we derive that
$$
\varphi^m_{m_1m'_1}\equiv {\sf H}_m (x_1^{m_1}x_{1'}^{m'_1})=\sum
_{k=0}^{{\rm min}(m_1,m'_1)} C^k_{m_1m'_1} Q^k
x_1^{m_1-k}x_{1'}^{m'_1-k}, \eqno (35)
$$
where
$$
 C^k_{m_1m'_1}=\frac{(q^{-2m_1};q^2)_k(q^{-2m'_1};q^2)_k}{(q^2;q^2)_k
(q^{-N-2m+4};q^2)_k} \frac{q^{(-n+2+\epsilon)k}}{(1+q^{N-2})^{k}}.
 \eqno (36)
$$

Using formula (16) the polynomials $\varphi^m_{m_1m'_1}$ can be
represented as
$$
\varphi^m_{m_1m'_1}=x_1^{m_1-m'_1}\sum_{k=0}^{m'_1} C^k_{m_1m'_1}
Q^k {Q'}^k(Q'_2/Q';q^2)_{m'-k}  \eqno (35')
$$
if $m_1\ge m'_1$ and as
$$
\varphi^m_{m_1m'_1}=\left( \sum_{k=0}^{m_1} C^k_{m_1m'_1} Q^k
{Q'}^k(Q'_2/Q';q^2)_{m'-k}\right) x_{1'}^{m'_1-m_1} \eqno (35'')
$$
if $m'_1\ge m_1$, where $Q'\equiv Q'_1$ and $Q'_2$ are such as in
(16).
\medskip

{\bf Theorem 1.} {\it The zonal polynomials $\varphi^m_{m_1m'_1}$,
$m_1+m'_1=m$, of ${\cal H}_m$ are orthogonal with respect to the
scalar product introduced in section VI. These polynomials
constitute a full set of zonal polynomials in the space ${\cal
H}_m$}
\medskip

{\it Proof.} We have ${\hat K}_1\triangleright
(x_1^{m_1}x_{1'}^{m'_1})=q^{m_1-m'_1} (x_1^{m_1}x_{1'}^{m'_1})$,
that is, the monomials $x_1^{m_1}x_{1'}^{m'_1}$, $m_1+m'_1=m$, are
eigenfunctions of the operator defined by the action of ${\hat
K}_1$ on ${\cal A}$ which belong to different eigenvalues. Since
the projection ${\sf H}_{m} :{\cal A}_{m} \to {\cal H}_{m}$
commutes with the action of $U_q({\rm so}_N)$, then ${\hat
K}_1\triangleright
\varphi^m_{m_1m'_1}=q^{m_1-m'_1}\varphi^m_{m_1m'_1}$. The scalar
product of section VI is defined in terms of the invariant
functional, that is, this scalar product is invariant with respect
to the action of ${\hat K}_i$, $i=1,2,\cdots ,n$. Since the zonal
polynomials $\varphi^m_{m_1m'_1}$, $m_1+m'_1=m$, belong to
different eigenvalues of ${\hat K}_1$, they are orthogonal.
Theorem is proved.
\medskip

It is possible to define zonal polynomials of the space ${\cal
H}_m$ with respect to the subalgebra $A:=U_q({\rm so}_2)\times
U_q({\rm so}_{N-2})$, where $U_q({\rm so}_2)$ is the subalgebra of
$U_q({\rm so}_N)$ generated by the element ${\hat K}_1$. Then the
following assertions are true.
\medskip

{\bf Theorem 2.} {\it The subspace of zonal polynomials of the
space ${\cal H}_m$ with respect to the subalgebra $A$ is not more
than one-dimensional. The space ${\cal H}_m$ contains a zonal
polynomial if and only if $m$ is even. This zonal polynomial
coincides with the polynomial $\varphi^m_{m/2,m/2}$ given by
formula (35).}
\medskip

Proof easily follows from the above results.
\bigskip

\noindent {\sf IX. ASSOCIATED SPHERICAL POLYNOMIALS WITH RESPECT
TO $SO_q(N-2)$}
 \medskip

The aim of this section is to construct an orthogonal basis of the
space ${\cal H}_{m}$ of homogeneous $q$-harmonic polynomials which
corresponds to the chain
$$
U_q({\rm so}_N)\supset U_q({\rm so}_{N-2})\supset \cdots \supset
U_q({\rm so}_3)\ \ ({\rm or}\ \ U_q({\rm so}_2)) .
$$
This basis is a $q$-analogue of the set of associated spherical
harmonics on the classical Euclidean space which are products of
Jacobi polynomials and correspond to the chain of the subgroups
$SO(N)\supset SO(2)\times SO(N-2)\supset \cdots$ (see [1], Chap.
10). The basis elements give solutions of the equation $\Delta
_qp=0$ in "separated coordinates". So, we obtain a $q$-analogue of
the classical separation of variables.

Let us note that
$$
\Delta_q\equiv \Delta_q^{(N)}=\sum_{j=1}^N q^{\rho_{j}}\partial_j
\partial_{j'}= (q^{\rho_{1}}\partial_1
\partial_{1'}+q^{\rho_{1'}}\partial_{1'}
\partial_{1})+ \Delta_q^{(N-2)},  \eqno (37)
$$
where $\Delta_q^{(N-2)}$ is the $q$-Laplace operator on the
subspace ${\cal A}^{(N-2)}\equiv {\Bbb C}_q[x_2,\cdots ,x_{2'}]$.
We also have from (9) that
$$
\partial_{1'} \partial_{1}-\partial_{1} \partial_{1'}=-\frac{q-q^{-1}}
{q^{\rho_{1}-1}+q^{-\rho_{1}+1}}\Delta_q^{(N-2)}. \eqno (38)
$$

Let $p(x_2,\ldots,x_{2'})$ be a polynomial of ${\cal A}$ which
does not depend on $x_1$ and $x_{1'}\equiv x_N$. Then it is easy
to see from (6) that $\partial_1 p(x_2,\ldots,x_{2'})=0$.
\medskip

{\bf Lemma 1.} {\it Let $p(x_2,\ldots,x_{2'})\in {\cal A}$ and
$\Delta_q^{(N-2)} p=0$. Then $\partial_{1'} p=0$ and $\Delta_q
p=0$.}
\medskip

{\it Proof.} Let $p(x_2,\ldots,x_{2'})$ is harmonic with respect
to $\Delta_q^{(N-2)}$, that is $\Delta_q^{(N-2)} p=0$. Then due to
(38) we have $\partial_{1} \partial_{1'} p=0$, and $\Delta_q p=0$
by (37). From formula (7) for $\partial_{1'}$, it follows that
$\partial_{1'} p=x_1 p'(x_2,\ldots,x_{2'})$, where
$p'(x_2,\ldots,x_{2'})$ is a polynomial in
$x_2,x_3,\ldots,x_{2'}$. Let us show that from $\partial_{1}
\partial_{1'} p=0$ the equality $\partial_{1'} p=0$ follows.
Indeed, due to (6) we have $0=\partial_{1}
\partial_{1'} p=\partial_{1} x_1 p'= \tilde{p}(x_2,\ldots,x_{2'})$,
where
$\tilde{p}(x_2,\ldots,x_{2'})$ is some polynomial in
$x_2,x_3,\ldots,x_{2'}$ which is a linear combination (with
nonvanishing coefficients) of the same monomials as $p'$ does.
Moreover, if $p'\ne 0$ then $\tilde{p}\ne 0$. Since $\tilde{p}= 0$
then $p'= 0$ and we have $\partial_{1'} p=x_1
p'(x_2,\ldots,x_{2'})=0$. It proves the lemma.
 \medskip

{\bf Lemma 2.} {\it If $p(x_2,\ldots,x_{2'})\in {\cal A}$, then
$$
\Delta_q^{(N-2)}(x_1^{m_1}x_{1'}^{m'_1}p(x_2,\ldots,x_{2'}))
=x_1^{m_1}x_{1'}^{m'_1}\Delta_q^{(N-2)}p(x_2,\ldots,x_{2'}). \eqno
(39)
$$
If $p(x_2,\ldots,x_{2'})$ is $\Delta_q^{(N-2)}$-harmonic, then}
$$
\Delta_q(x_1^{m_1}x_{1'}^{m'_1}p(x_2,\ldots,x_{2'}))= (q^{\rho_1}+
q^{-\rho_1})\partial_{1'} \partial_{1}
x_1^{m_1}x_{1'}^{m'_1}p(x_2,\ldots,x_{2'}). \eqno (40)
$$

{\it Proof.} Since $\partial_1 \hat x_{1'}=q \hat x_{1'}
\partial_1$, $\partial_2 \hat x_{1}= \hat x_{1}
\partial_2$, $\partial_1 p(x_2,\ldots,x_{2'})=0$, and
$$
\partial_2{\hat x}_{1'}=
{\hat x}_{1'}\partial_2+(q-q^{-1})q^{\rho_1-\rho_2}
{\hat x}_{2'}\partial_1
$$
(see formulas (11)-(13)) we obtain
$$
\partial_2 (x_1^{m_1}x_{1'}^{m'_1} p(x_2,\cdots,x_{2'}))=
\hat x_1^{m_1} \partial_2 (x_{1'}^{m'_1} p(x_2,\ldots,x_{2'}))
$$  $$
= \hat x_1^{m_1} (\hat x_{1'}
\partial_2+(q-q^{-1})q^{\rho_1-\rho_{2}}\hat x_{2'} \partial_1)
x_{1'}^{m'_1-1} p(x_2,\cdots,x_{2'})
$$  $$
= \hat x_1^{m_1}\hat
x_{1'}^{m'_1} \partial_2 p(x_2,\cdots,x_{2'}).
$$
Analogously, using relation (12) we derive that
$$
\partial_{2'} (x_1^{m_1}x_{1'}^{m'_1} p(x_2,\ldots,x_{2'}))=
\hat x_1^{m_1} \partial_{2'} (x_{1'}^{m'_1} p(x_2,\ldots,x_{2'}))
$$  $$
= \hat x_1^{m_1} (\hat x_{1'}
\partial_{2'}+(q-q^{-1})q^{\rho_1-\rho_{2'}}\hat x_{2} \partial_1)
(x_{1'}^{m'_1-1} p(x_2,\ldots,x_{2'}))
$$  $$
= \hat x_1^{m_1}\hat x_{1'}^{m'_1} \partial_{2'}
p(x_2,\cdots,x_{2'}).
$$
We have the same results when $\partial_2$ and $\partial_{2'}$ are
replaced by $\partial_i$ and $\partial_{i'}$, $i=3,4,\cdots$. This
leads to the relation (39). If $p$ is $\Delta_q^{(N-2)}$-harmonic,
then it follows from (38) that
$(\partial_{1'}\partial_1-\partial_1
\partial_{1'}) x_1^{m_1}x_{1'}^{m'_1} p(x_2,\cdots,x_{2'})=0$.
From here and from (37) we derive that
$$
\Delta_q( x_1^{m_1}x_{1'}^{m'_1} p)=(
q^{\rho_1}\partial_{1}\partial_{1'}+
q^{\rho_{1'}}\partial_{1'}\partial_{1})
( x_1^{m_1}x_{1'}^{m'_1} p)=
(q^{\rho_1}+q^{\rho_{1'}})\partial_{1'}\partial_{1}(
x_1^{m_1}x_{1'}^{m'_1} p)
$$
and the relation (40) is proved. Lemma is proved.
\medskip

{\bf Proposition 6.} {\it Let $h\equiv h(x_2,\cdots,x_{2'})$ be a
$\Delta_q^{(N-2)}$-harmonic polynomial of degree $l$. Then}
$$
\Delta_q (x_1^{m_1}x_{1'}^{m'_1}
h)=(q^{\rho_1}+q^{-\rho_1})[m_1]_q[m'_1]_q
q^{m_1+m'_1-1}x_1^{m_1-1}x_{1'}^{m'_1-1} h. \eqno (41)
$$

{\it Proof.} Using (40) and then (1), (6), (13), we derive that
$$
 \Delta_q (x_1^{m_1}x_{1'}^{m'_1} h(x_2,\ldots,x_{2'}))=
(q^{\rho_1}+q^{-\rho_1}) [m_1]_q q^{m'_1+m_1+l-1} x_1^{m_1-1}
\partial_{1'} x_{1'}^{m'_1}h(x_2,\ldots,x_{2'}). \eqno (42)
$$
By (11) we have for $\partial_{1'} x_{1'}^{m'_1}h$ the expression
$$
\partial_{1'} x_{1'}^{m'_1}h(x_2,\cdots,x_{2'})=
\Bigl( q^{-1} x_{1'} \partial_{1'}-(q-q^{-1}) \sum_{j<N} {\hat
x}_j\partial_j
$$  $$
 +(q-q^{-1}) {\hat x}_1\partial_1
q^{2\rho_1}+c\Bigr) x_{1'}^{m'_1-1}h(x_2,\ldots,x_{2'})
$$  $$
= \Bigl( q {\hat x}_{1'} \partial_{1'}-(q-q^{-1}) E
+q^{m'_1-1+l}\Bigr) x_{1'}^{m'_1-1}h(x_2,\cdots,x_{2'}), \eqno
(43)
$$
where $E=\sum _{k=1}^N {\hat x}_k\partial_k$. It is proved by
using the relation between ${\hat x}_i$ and $\partial_j$ that
$$
E{\hat x}_k=q^{-1}{\hat x}_kE+\frac{q-q^{-1}}{1+q^{N-2}}
q^{N-\rho_k-2}{\hat Q}\partial_{k'}+{\hat x}_kc
$$
(see also [3], Proposition 2.9). Then
$$
E(x_{1'}^{m'_1}h(x_2,\ldots,x_{2'}))=(q^{-1}{\hat
x}_{1'}E+x_{1'}q^{l+m'_1-1}) x_{1'}^{m'_1-1}h(x_2,\ldots,x_{2'})
$$   $$
= (q^{l+m'_1-1}+q^{l+m'_1-3}+\cdots
+q^{l-m'_1+1})x_{1'}^{m'_1}h(x_2,\ldots,x_{2'})+ q^{-m'_1} {\hat
x}_{1'}^{m'_1}E h(x_2,\ldots,x_{2'}).
$$
By direct calculation it is proved that
$$
E=\frac{c-c^{-1}}{q-q^{-1}}+\frac{q-q^{-1}}{(1+q^{N-2})^2}
q^{N-1}{\hat Q}\Delta_qc^{-1}
$$
(see also [3]). Since $h\in {\cal H}_l$, then $Eh=[l]_qh$. Now we
have for $E{\hat x}_{1'}^{m'_1}h$ the expression
$$
E(x_{1'}^{m'_1}h(x_2,\ldots,x_{2'}))=(q^l[m'_1]_q+q^{-m'_1}[l]_q)
x_{1'}^{m'_1}h(x_2,\ldots,x_{2'})= [m'_1+l]_q
x_{1'}^{m'_1}h(x_2,\ldots,x_{2'}).
$$
Therefore, returning to (43) we obtain
$$
\partial_{1'} x_{1'}^{m'_1}h(x_2,\cdots,x_{2'})=
\Bigl( q {\hat x}_{1'}
\partial_{1'}+q^{m'_1-1+l}-(q-q^{-1})[m'_1-1+l]_q\Bigr)
x_{1'}^{m'_1-1}h(x_2,\ldots,x_{2'})
 $$   $$
=(q {\hat x}_{1'} \partial_{1'}+q^{-m'_1-l+1})
x_{1'}^{m'_1-1}h(x_2,\ldots,x_{2'}).
$$
Applying this relation for  $x_{1'}^{m'_1-1}h$,
$x_{1'}^{m'_1-2}h,\cdots, x_{1'}h$ and Lemma 1 we receive
$$
\partial_{1'} x_{1'}^{m'_1}h(x_2,\cdots,x_{2'})
=(q^{1-m'_1-l}+q^{1-m'_1-l+2}+\cdots)x_{1'}^{m'_1-1}
h(x_2,\cdots,x_{2'})
$$  $$
+ q^{m'_1} x_{1'}^{m'_1}\partial_{1'}h(x_2,\cdots,x_{2'}) =
q^{-l}[m'_1]_q x_{1'}^{m'_1-1}h(x_2,\ldots,x_{2'}).
$$
Now using (42) we derive (41). Proposition is proved.
 \medskip

Let $h_l\in {\cal H}_l^{(N-2)}$, where ${\cal H}_l^{(N-2)}$ is the
space of $\Delta_q^{(N-2)}$-harmonic polynomials in
$x_2,x_3,\cdots ,x_{2'}$. Then $x_{1}^{m_1}x_{1'}^{m'_1} h_l \in
{\cal A}_m$, $m=m_1+m'_1+l$. Using formula (41), in the same way
as in the case of formula (34), we find that
$$
\Delta^k_q(x_1^{m_1}x_{1'}^{m'_1}h_l)=(1+q^{N-2})^kq^{(m_1+m'_1-k)k}
q^{-(n-\epsilon)k} \frac{[m_1]_q![m'_1]_q!}{[m_1-k]_q![m'_1-k]_q!}
x_1^{m_1-k}x_{1'}^{m'_1-k}h_l,  \eqno (44)
$$
where $\epsilon=1$ for $N=2n$ and $\epsilon =\frac 12$ for
$N=2n+1$. Now using formulas (30) and (32) we derive that
$$
{\sf H}_m (x_1^{m_1}x_{1'}^{m'_1}h_l)=t^{N,m}_{m_1m'_1}h_l, \eqno
(45)
$$
where
$$
t^{N,m}_{m_1m'_1}= \sum _{k=0}^{{\rm min}(m_1,m'_1)}
C^{m,k}_{m_1m'_1} Q^k x_1^{m_1-k}x_{1'}^{m'_1-k}, \eqno (46)
$$  $$
 C^{m,k}_{m_1m'_1}=\frac{(q^{-2m_1};q^2)_k(q^{-2m'_1};q^2)_k}
{(q^2;q^2)_k (q^{-N-2m+4};q^2)_k}
\frac{q^{(-n+\epsilon-2l+2)k}}{(1+q^{N-2})^{k}} .
$$
By means of formula (16) the polynomials $t^{N,m}_{m_1m'_1}$ can
be represented in the form similar to (35$'$) and (35$''$).
\medskip

{\bf Proposition 7.} {\it The space ${\cal H}_m$ can be
represented as the orthogonal sum
$$
{\cal H}_m=\bigoplus _{m_1,m'_1} t^{N,m}_{m_1m'_1}{\cal
H}^{(N-2)}_{m-m_1-m'_1},  \eqno (47)
$$
where ${\cal H}^{(N-2)}_{m-m_1-m'_1}$ is the space of
$\Delta^{(N-2)}_q$-harmonic polynomials in $x_2,x_3,\cdots
,x_{2'}$ and summation is over all nonnegative values of $m_1$ and
$m'_1$ such that $m-m_1-m'_1\ge 0$.}
 \medskip

{\it Proof.} The subspaces $t^{N,m}_{m_1m'_1}{\cal
H}^{(N-2)}_{m-m_1-m'_1}$ from (47) do not pairwise intersect and
elements from different subspaces are linearly independent.
Therefore, on the right hand side of (47) we have a direct sum.
Besides, we have ${\cal H}_m\supseteq \bigoplus _{m_1,m'_1}
t^{N,m}_{m_1m'_1}{\cal H}^{(N-2)}_{m-m_1-m'_1}$. By direct
computation (by using Corollary 3) we show that dimensions of
spaces on both sides of (47) are equal to each other. Now in order
to prove our proposition we have to show that the sum on the right
hand side is orthogonal.

It is easy to prove that the subspaces on the right hand side of
(47) are eigenspaces of operators ${\hat K}_i$, $i\le n$, from
formula (5) belonging to different eigenvalues. As in the proof of
Theorem 1 it follows from here that the sum (47) is orthogonal.
Proposition is proved.
\medskip

Now we apply the decomposition (47) to the subspaces ${\cal
H}^{(N-2)}_{m-m_1-m'_1}$ and obtain
$$
{\cal H}_m= \bigoplus _{m_1,m'_1} \bigoplus _{m_2,m'_2}
t^{N,m}_{m_1m'_1} t^{N-2,l}_{m_2m'_2}{\cal
H}^{(N-4)}_{l-m_2-m'_2}, \ \ \ l=m-m_1-m'_1,
$$
where ${\cal H}^{(N-4)}_{l-m_2-m'_2}$ are the subspaces of
homogeneous $q$-harmonic polynomials in $x_3,x_4,\cdots ,x_{3'}$.
Continuing such decompositions we obtain the decomposition
$$
{\cal H}_m= \bigoplus _{{\bf m},{\bf m}',k} {\Bbb C}\Xi _{{\bf
m},{\bf m}',k}(x_1,\cdots ,x_{1'})   \eqno (48)
$$
if $N=2n$ and the decomposition
$$
{\cal H}_m= \bigoplus _{{\bf m},{\bf m}',\sigma} {\Bbb C}\Xi
_{{\bf m},{\bf m}',\sigma}(x_1,\cdots ,x_{1'})   \eqno (49)
$$
if $N=2n+1$, where ${\bf m}=(m_1,m_2,\cdots ,m_{n-1})$, ${\bf
m}'=(m'_1,m'_2,\cdots ,m'_{n-1})$ in the first case, ${\bf
m}=(m_1,m_2,\cdots ,m_{n})$, ${\bf m}'=(m'_1,m'_2,\cdots ,m'_{n})$
in the second case and $m_{j}$ are nonnegative integers, $k$ take
integral values, and $\sigma=0$ or 1. The basis $q$-harmonic
polynomials $\Xi _{{\bf m},{\bf m}',k}(x_1,\cdots ,x_{1'})$ of
${\cal H}_m$ are given by the formula
$$
\Xi _{{\bf m},{\bf m}',k}(x_1,\cdots ,x_{1'})= t^{N,m}_{m_1m'_1}
t^{N-2,m-m_1-m'_1}_{m_2m'_2}\cdots t^{4,m-\sum_{i=1}^{n-2}m_i
-\sum_{i=1}^{n-2}m'_i}_{m_{n-1}m_{n-1}'}t^{2,k}  \eqno (50)
$$
if $N=2n$ and by the formula
$$
\Xi _{{\bf m},{\bf m}',\sigma }(x_1,\cdots ,x_{1'})=
t^{N,m}_{m_1m'_1} t^{N-2,m-m_1-m'_1}_{m_2m'_2}\cdots
t^{3,m-\sum_{i=1}^{n-1}m_i
-\sum_{i=1}^{n-1}m'_i}_{m_{n}m_{n}'}x_{n+1}^\sigma   \eqno (51)
$$
if $N=2n+1$. In (50), $t^{2,k}=x_n^k$ if $k>0$, $t^{2,k}=1$ if
$k=0$, and $t^{2,k}=x_{n'}^{-k}$ if $k<0$. Note that the integers
$k$, $\sigma$, ${\bf m}=(m_1,m_2,\cdots )$ and ${\bf
m}'=(m'_1,m'_2,\cdots )$ take such values that
$$
m_1+m'_1+m_2+m'_2+\cdots +m_{n-1}+m'_{n-1}+k=m   \eqno (52)
$$
for $N=2n$ and
$$
m_1+m'_1+m_2+m'_2+\cdots +m_{n}+m'_{n}+\sigma =m   \eqno (53)
$$
for $N=2n+1$. Besides, conditions like the condition
$m-m_1-m'_1\ge 0$ of Proposition 7 must be fulfilled on each step.
\medskip

{\bf Theorem 3.} {\it If $N=2n$ then the polynomials (50) for
which the equality (52) is satisfied constitute an orthogonal
basis of the space ${\cal H}_m$. If $N=2n+1$ then the polynomials
(51) for which the equality (53) is satisfied constitute an
orthogonal basis of the space ${\cal H}_m$.}
\medskip

{\it Proof.} The fact that the polynomials (50) for $N=2n$ and the
polynomials (51) for $N=2n+1$ constitute a basis of ${\cal H}_m$
was proved above. Orthogonality of basis elements is proved in the
same method as in Theorem 1. Theorem is proved.
\medskip

It is interesting to have explicit formulas how the generators
$K_i, E_i, F_i$ of $U_q({\rm so}_n)$ act on the basis elements of
Theorem 3. However, derivation of these formulas are very awkward
and the formulas are not simple. We shall consider them in a
separate paper.


\begin{thebibliography}{99}

\bibitem{1}
Vilenkin, N.~Ja. and Klimyk, A.~U. {\it Representation of Lie
Groups and Special Functions}, vol.~2 (Dordrecht: Kluwer Academic
Publishers), 1993.

\bibitem{2}Iorgov, N. Z. and Klimyk, A. U., A Laplace operator and
harmonics on the quantum complex vector space, {\it J. Math.
Phys.} {\bf 44} (2003), 823--848.

\bibitem{3}Sugitani T., Harmonic analysis on quantum spheres
associated with representations of $U_q({\rm so}_N)$ and
$q$-Jacobi polynomials, {\it Compos. Math.} {\bf 99} (1995),
249--281.

\bibitem{4}Reshetikhin, N. Ya., Takhtajan, L. A., and Faddeev, L. D.,
Quantization of Lie groups and Lie algebras, {\it Leningrad Math.
J.} {\bf 1} (1990), 193-225.

\bibitem{5} Klimyk, A. and Schm\"udgen, K., {\it Quantum Groups and
Their Representations} (Berlin: Springer), 1997.

\bibitem{6} Kostant, B., Lie group
representations on polynomial rings,
{\it Amer. J. Math.} {\bf 83} (1963), 327--404.

\bibitem{7} Burban, I. M. and Klimyk, A. U., Representations of
the quantum
algebra $U_q({\rm su}_{1,1})$, {\it J. Phys. A: Math. Gen.} {\bf
26} (1993), 2139--2151.

\bibitem{8} Noumi, M., Macdonald's symmetric polynomials as zonal
spherical functions on some quantum homogeneous spaces, {\it Adv.
Math.} {\bf 123} (1996), 16--77.

\bibitem{9} Gasper, G. and Rahman, M., {\it Basic Hypergeometric
Functions} (Cambridge: Cambridge Univ. Press), 1990.

\end{thebibliography}
\end{document}